\input amssym.def
\input epsf

\let \blskip = \baselineskip
\parskip=1.2ex plus .2ex minus .1ex

\tabskip 20pt
\tolerance = 1000
\pretolerance = 50
\newcount\itemnum
\itemnum = 0
\overfullrule = 0pt

\def\title#1{\bigskip\centerline{\bigbigbf#1}}
\def\author#1{\bigskip\centerline{\bf #1}\smallskip}
\def\address#1{\centerline{\it#1}}
\def\abstract#1{\vskip1truecm{\narrower\noindent{\bf Abstract.} #1\bigskip}}

\def\sp{\bigskip}
\def\nosp{\vskip -\the\blskip plus 1pt minus 1pt}

\def\br{\hfil\break} 
\def\ti{\br \hglue \the \parindent}

\def\skipit#1{}
\def\mdag{\raise 3pt\hbox{\dag}}

\def\IP{\par\hang}
\def\XP{\par\noindent\hang}
\def\LP{\par\noindent}
\def\BP[#1]{\par\item{[#1]}}
\def\SH#1{\sp\vskip\parskip\leftline{\bigbf #1}\nobreak}

\def\TH#1{\sp\XP{\bf THEOREM\ \shead#1}}

\def\CO#1{\sp\XP{\bf COROLLARY\ \shead#1}}

\def\PF{\LP{\bf Proof:\ }}

\def\NX{\advance\itemnum by 1 \sp\LP {\bf \shead \the\itemnum.\ }}
\def\qed{\null\nobreak\hfill\hbox{${\vrule width 5pt height 6pt}$}\par\sp}

\def\cart{\>\hbox{${\vcenter{\vbox{
    \hrule height 0.4pt\hbox{\vrule width 0.4pt height 4.5pt
    \kern4pt\vrule width 0.4pt}\hrule height 0.4pt}}}$}\>}
\def\bxmu{\>\hbox{${\vcenter{\vbox {
    \hrule height 0.4pt\hbox{\vrule width 0.4pt height 4pt
    \hskip -1.3pt\lower 1.8pt\hbox{$\times$}\negthinspace\vrule width 0.4pt}
    \hrule height 0.4pt}}}$}\>}

\def\lin#1{\hbox to #1true in{\hrulefill}}



\def\al{\alpha}




\def\({\left(}	\def\){\right)}



\def\SET#1:#2{\{#1\colon\;#2\}}

		
	\def\un#1{\underline{#1}}



\magnification=\magstep1
\vsize=9.0 true in
\hsize=6.5 true in
\headline={\hfil\ifnum\pageno=1\else\folio\fi\hfil}
\footline={\hfil\ifnum\pageno=1\folio\else\fi\hfil}

\parindent=20pt
\baselineskip=12pt
\parskip=.5ex  

\def\shead{ }

\font\bigbf = cmb10 scaled \magstep1

\font\bigbigbf = cmb10 scaled \magstep2




\title{A Short Proof that ``Proper $=$ Unit''}
\author{Kenneth P. Bogart}
\address{Dartmouth College, Hanover, NH 03755-3551, {k.p.bogart@dartmouth.edu}}
\author{Douglas B. West\mdag}
\address{University of Illinois, Urbana, IL 61801-2975, {west@math.uiuc.edu}}
\vfootnote{}{\br
   \mdag Research supported in part by NSA/MSP Grant MDA904-93-H-3040.\br
   Running head: Proper $=$ Unit \br
}
\abstract{A short proof is given that the graphs with proper interval
representations are the same as the graphs with unit interval representations.}

\sp
An graph is an {\it interval graph} if its vertices can be assigned intervals
on the real line so that vertices are adjacent if and only if the corresponding
intervals intersect; such an assignment is an {\it interval representation}.
When the intervals have the same length, we have a {\it unit interval
representation}.  When no interval properly contains another, we have a
{\it proper interval representation}.  The {\it unit interval graphs} and
{\it proper interval graphs} are the interval graphs having unit interval
or proper interval representations, respectively.

Since no interval contains another of the same length, every unit interval
graph is a proper interval graph.  Roberts [1] proved that also every proper
interval graph is a unit interval graph; the two classes are the same.
He proved this as part of a characterization of unit interval graphs as the
interval graphs with no induced subgraph isomorphic to the ``claw'' $K_{1,3}$
(it is immediate that the condition is necessary).  Roberts used a version of
the Scott-Suppes [2] characterization of semiorders to prove that claw-free
interval graphs are unit interval graphs.  By eschewing the trivial implication
that unit interval graphs are proper interval graphs and instead going from
``claw-free'' to ``proper'' to ``unit'' among the interval graphs,
we obtain a short self-contained proof.

In the language of partial orders, our proof also characterizes the semiorders
among the interval orders.  A partial order is an {\it interval order} if its
elements can be assigned intervals on the real line so that $x<y$ if and only if
the interval assigned to $x$ is completely to the left of the interval assigned
to $y$.  A partial order is a {\it semiorder} if its elements can be assigned
numbers so that $x<y$ if and only if the number assigned to $y$ exceeds the
number assigned to $x$ by more than 1.  The poset $\un1+\un3$ is the poset
consisting of two disjoint chains of sizes 3 and 1.  The semiorders are
precisely the interval orders that do not contain $\un1+\un3$.

\vfil
\eject
\TH{}
The following statements are equivalent when $G$ is a simple graph.
\IP
A) $G$ is a unit interval graph.
\br 
B) $G$ is an interval graph with no induced $K_{1,3}$.
\br 
C) $G$ is a proper interval graph.
\PF
In an interval representation of $K_{1,3}$, the intervals for the three leaves
must be pairwise disjoint, and then the interval for the central vertex must
properly contain the middle of the three intervals for leaves.  Thus A and
C imply B.

For the converse, let $G$ be a claw-free interval graph, and consider an
interval representation that assigns to each $v\in V(G)$ an interval $I_v$.
We first transform this into a proper interval representation.  Since $G$ is
claw-free, there is no pair $x,y\in V(G)$ such that 1) $I_y\subset I_x$ and 2)
$I_x$ intersects intervals to the left and right of $I_y$ that do not intersect
$I_y$.  If $I_x=[a,b]$ and $I_y=[c,d]$ with $a<c\le d<b$, this means that
$[a,c]$ or $[d,b]$ is empty of endpoints of intervals that don't intersect
$I_y$.  Hence we can extend $I_y$ past the end of $I_x$ on one end without
changing the graph obtained from the representation.  Repeating this until no
more pairs of intervals are related by inclusion yields a proper interval
representation.

From a proper interval representation of $G$, we obtain a unit interval
representation.  When no interval properly includes another, the right endpoints
have the same order as the left endpoints.  We process the representation from
left to right, adjusting all intervals to length 1.  At each step until all
have been adjusted, let $I_x=[a,b]$ be the unadjusted interval that has the
leftmost left endpoint.  Let $\al=a$ unless $I_x$ contains the right endpoint of
some other interval, in which case let $\al$ be the largest such right endpoint.
Such an endpoint would belong to an interval that has already been adjusted to
have length 1; thus $\al<\min\{a+1,b\}$.  Now, adjust the portion of the
representation in $[a,\infty)$ by shrinking or expanding $[\al ,b]$ to
$[\al ,a+1]$ and translating $[b,\infty)$ to $[a+1,\infty)$.  The order of
endpoints does not change, intervals earlier than $I_x$ still have length 1, and
$I_x$ also now has length 1.  Iterating this operation produces the unit
interval representation.  \qed

This theorem has a standard interpretation for posets.
The incomparability graph of an interval order is an interval graph.
Existence of the function required in the definition of semiorder
is equivalent to having a representation as an interval order using
intervals of unit length.  The incomparability graph of $\un1+\un3$ is
$K_{1,3}$.  Thus applying the Theorem above to incomparability graphs 
yields the following Corollary.

\CO{}
The following statements are equivalent when $P$ is a poset.
\br
A) $P$ is a unit interval order.
\br
B) $P$ is a semiorder.
\br
C) $P$ is a interval order not containing $\un1+\un3$.
\br
D) $P$ is a proper interval order.

\SH
{References}
\frenchspacing
\BP [1]
F.S.~Roberts, Indifference graphs, {\it Proof Techniques in Graph Theory}
(F.~Harary, ed.).  Academic Press (1969), 139--146.
\BP [2]
D.~Scott and P.~Suppes, Foundational aspects of theories of measurement,
{\it J. Symbolic Logic} 23(1958), 233--247.

\bye